\documentclass[12pt,A4,reqno]{amsart}
\usepackage{amsfonts, amssymb, amsmath, amsthm, }
\theoremstyle{plain}

    \theoremstyle{definition}

\newcommand{\R}{\mathbb R}
\newcommand{\Z}{\mathbb Z}

\newcommand{\n}{\vskip .1cm}
\newcommand{\nn}{\vskip .2cm}

\begin{document}
\title{On the equivalence of several definitions of compact infra-solvmanifolds}

\author{Shintar\^{o} Kuroki$^\dag$}
\address{$^\dag$ Osaka City University Advanced Mathematical Institute (OCAMI),
      3-3-138 Sugimoto, Sumiyoshi-ku, Osaka 558-8585, Japan \newline \textit{and}
             Department of Mathematics, University of Toronto, Room 6290, 40 St.
     George Street, Toronto, Ontario, M5S 2E4, Canada}
 \email{kuroki@scisv.sci.osaka-cu.ac.jp, shintaro.kuroki@utoronto.ca}

\author{Li Yu*}
\address{*Department of Mathematics and IMS, Nanjing University, Nanjing, 210093, P.R.China}
 \email{yuli@nju.edu.cn}
 \thanks{2010 \textit{Mathematics Subject Classification}.  22E25, 22E40, 22F30, 53C30\\
 The first author was supported in part by the JSPS Strategic Young Researcher
   Overseas Visits Program for Accelerating Brain Circulation ``Deepening and Evolution of
     Mathematics and Physics, Building of International Network Hub based on OCAMI''.
 The second author is partially supported by
 Natural Science Foundation of China (Grant no.11001120) and the
 PAPD (priority academic program development) of Jiangsu higher education institutions.}

 \keywords{Infra-solvmanifold, solvable Lie group, homogeneous space,
  discrete group action, holonomy.
 }

 \begin{abstract}
   We show the equivalence of several definitions of
   compact infra-solvmanifolds that appear in various math literatures.
 \end{abstract}


\maketitle

 Infra-solvmanifolds are a speical class of aspherical manifolds
  studied by mathematicians for a long time.
 By name, an infra-solvmanifold is finitely covered by a solvmanifold
 (which is the quotient of a connected solvable Lie group by a closed subgroup).
 Infra-solvmanifolds include all
 flat Riemannian manifolds and infra-nilmanifolds.\nn

 Compact Infra-solvmanifolds are \emph{smoothly rigid}, i.e.
 any homotopy equivalence between two compact infra-solvmanifolds is
 homotopic to a diffeomorphism (see~\cite{FarJone97} and~\cite{Wilk00}).
 Geometrically, a compact infra-solvmanifold $M$ can be described (see~\cite[Proposition 3.1]{Tusch97})
 as a manifold admitting a sequence of
  Riemannian metric $\{ g_{n}\}$ with uniformly bounded sectional curvature
  so that $(M, g_n)$ collapses in the \emph{Gromov-Hausdorff} sense to a \emph{flat orbifold}.
 \nn

 In the history of study on this topic, there are several definitions of compact
 infra-solvmanifolds
 appearing in various math literatures. This is mainly because that people look at
 these manifolds from very different angles (topologically, algebraically or
 geometrically).\nn

  The following are three of these definitions.\nn

 \begin{itemize}
   \item[Def 1:] Let $G$ be a connected, simply connected solvable Lie group, $K$ be a maximal compact
 subgroup of the group $\mathrm{Aut}(G)$ of automorphisms of $G$, and $\Gamma$ be a
  cocompact, discrete subgroup of $G \rtimes K$. If the
  action of $\Gamma$ on $G$ is free and $[\Gamma: G\cap \Gamma ] < \infty$,
  the orbit space $\Gamma \backslash G$ is called an
  \emph{infra-solvmanifold modeled on $G$}.  See~\cite[Definition 1.1]{Tusch97}.\nn

   \item[Def 2:] A \emph{compact infra-solvmanifold} is a manifold of the form
    $\Gamma \backslash G$, where $G$ is a connected,
    simply connected solvable Lie group, and $\Gamma$ is a torsion-free
  cocompact discrete subgroup of $\mathrm{Aff}(G)=G \rtimes \mathrm{Aut}(G)$
   which satisfies: the closure of $hol(\Gamma)$ in $\mathrm{Aut}(G)$ is compact where
   $hol: \mathrm{Aff}(G) \rightarrow \mathrm{Aut}(G)$ is the holonomy
   projection.\nn

   \item[Def 3:] A \emph{compact infra-solvmanifold} is a double coset space $\Gamma \backslash G \slash K$
   where $G$ is a virtually connected and virtually solvable Lie
   group, $K$ is a maximal compact subgroup of $G$ and $\Gamma$ is a
   torsion-free, cocompact, discrete subgroup of $G$. See~\cite[Definition 2.10]{FarJone90}.
 \end{itemize}

 A \emph{virtually connected} Lie group is a Lie group with finitely many
 connected components. \nn

 The main purpose of this note is to explain why the above three definitions of
 compact infra-solvmanifolds are equivalent. The authors are fully aware that the
 reason is probably known to many people. But since we did not find any formal proof of
 this equivalence and felt this phenomenon a little confusing,
   so we want to write a proof here for the convenience of future
   reference. However, we do not intend to give a completely new treatment of this subject.
   Our proof will directly quote some results on infra-solvmanifolds from~\cite{Baube04} and~\cite{Wilk00}
  and use many subtle facts in Lie theory. This note can be
 treated as an elementary exposition of compact infra-solvmanifolds and some related concepts.
 \\

 \noindent $\textbf{\S 1.}$ \textbf{The equivalence of Def 1 and Def 2}\nn

  Let $M$ be a compact infra-solvmanifold in the sense of Def 1.
   First of all, any $g \in \Gamma $ can be decomposed as $g = k_g u_g $ where
   $k_g \in K \subset \mathrm{Aut}(G)$ and $u_g\in G$. The holonomy projection $hol:  G \rtimes
   \mathrm{Aut}(G)\rightarrow \mathrm{Aut}(G) $ sends
   $g$ to $k_g$. Since $hol$ is a group homomorphism, its image
   $hol(\Gamma)$ is a subgroup of $K$. By assumption,
   $ | hol(\Gamma) | = [\Gamma : G\cap \Gamma]$ is
   finite, so $hol(\Gamma)$ is compact.
   In addition,
   since $G$ is a connected simply-connected solvable Lie group, so $G$ is
    diffeomorphic to an Euclidean space. Then
    by Smith fixed point theorem (\cite[Theorem I]{Smith41}),
     $\Gamma$ acting on $G$ freely implies that $\Gamma$ is torsion-free.
   So $M$ satisfies Def 2.
   \nn
   Conversely, if $M$ is a compact infra-solvmanifold in Def 2,
   then~\cite[Theorem 3 (a) $\Rightarrow$ (f)]{Wilk00} tells us that
   there exists a connected, simply-connected
   solvable Lie group $G'$ so that the $\Gamma$ (which defines $M$)
   can be thought of as discrete cocompact subgroup of $G'\rtimes F$ where $F$ is a
   finite subgroup of $\mathrm{Aut}(G')$.
   Moreover, there exists an equivariant diffeomorphism from $G$ to $G'$ with respect to
   the action of $\Gamma$ (by~\cite[Theorem 1 and Theorem 2]{Wilk00}). Hence
    $$M = \Gamma \backslash G \underset{\mathrm{diff.}}{\cong}
     \Gamma \backslash G'.$$
  Then $[\Gamma : \Gamma\cap G'] = | hol_{G'}(\Gamma) | \leq  | F | $ is
  finite.\nn

  It remains to show that the action of $\Gamma$ on
  $G'$ is free.
  Since $F$ is finite, we can choose an $F$-invariant Riemannian metric
  on $G'$ (in fact we only need $F$ to be compact).
  If the action of an element $g\in \Gamma$ on $G'$ has a fixed point, say
  $h_0$. Let $L_{h_0}$ be the left translation of $G'$ by $h_0$.
   Then $L^{-1}_{h_0} g L_{h_0}$ fixed the identity element, which
   implies $L^{-1}_{h_0} g L_{h_0} \in F$. So
   there exists a compact neighborhood $U$ of $h_0$ so that $g\cdot U = U$ and
  $g$ acts isometrically on $U$ with respect to the Riemannian metric just as the
  element $L^{-1}_{h_0} g L_{h_0} \in F$ acts around $e$.
  So $A=Dg : T_{h_0} U \rightarrow T_{h_0} U$
  is an orthogonal transformation.
  Then $A$ is conjugate in $O(n)$ to the block diagonal matrix of the form
     \[
           \begin{pmatrix}
             B(\theta_1) &  \ & 0 \\
             \ &  \ddots & \ \\
             0 &  \ &  B(\theta_m)
           \end{pmatrix}
     \]
   where $B(0)=1$, $B(\pi)=-1$, and $B(\theta_j)=
\left(\begin{array}{cc}
          \cos \theta_j & -\sin \theta_j\\
          \sin \theta_j & \cos \theta_j
         \end{array}\right)$, $1\leq j \leq m$.\nn

  Since $\Gamma$ is torsion-free, $g$ is an infinite order element.
  This implies that at least one $\theta_{j}$ is irrational. Then
   for any $v\neq 0 \in T_{h_0}G'$, $| \{ A^n v \}_{n\in \Z} | =\infty$.
   So there exists $h\in U$ so that $|\{ g^n \cdot h \}_{n\in \Z} | = \infty$. Then the
    set $\{ g^n \cdot h \}_{n\in \Z} \subset U$ has
   at least one accumulation point. This contradicts the fact that
  the orbit space $\Gamma \backslash G'$ is Hausdorff (since
   $\Gamma\backslash G' = M$ is a manifold).\nn

  Combing the above arguments, $M$ is an infra-solvmanifold modeled on $G'$ in the sense of Def 1.
  \hfill $\qed$\nn

   In the following, by writing ``$\mathbf{A} \Rightarrow \mathbf{B}$'' we mean
   a space satisfying condition $\mathbf{A}$ implies that it satisfies condition $\mathbf{B}$.
    \\

 \noindent $\textbf{\S 2.}$ \textbf{Def 1 $\Rightarrow$ Def 3} \nn

   Let $M$ be an infra-solvmanifold in the sense of Def 1.
   Let $hol(\Gamma)$ be the image of
   the the holonomy projection $hol: \Gamma \rightarrow \mathrm{Aut}(G)$. Then $hol(\Gamma)$ is a
   finite subgroup of $\mathrm{Aut}(G)$.\nn

    Define $\widetilde{G} = G \rtimes hol(\Gamma)$ which is virtually solvable.
    It is easy to see that $hol(\Gamma)$ is
  a maximal compact subgroup of $\widetilde{G}$ and, $\Gamma$ is a
  cocompact, discrete subgroup of $\widetilde{G}$.
    Then $M \cong \Gamma\backslash \widetilde{G} \slash
   K$. So $M$ satisfies Def 3. \hfill $\qed$ \\

\noindent $\textbf{\S 3.}$ \textbf{Def 3 $\Rightarrow$ Def 2}
 \nn
   Let $M=\Gamma \backslash G \slash K$ be a compact infra-solvmanifold in the sense of Def 3.
   Since $K$ is a maximal compact subgroup of $G$, so $G\slash K$ is
   contractible. Note that $K$ is not necessarily a normal subgroup of $G$, so $G\slash K$ may not directly
   inherit a group structure from $G$.
   \nn

    Let $G_0$ be the connected component of $G$ containing the
   identity element. Then by~\cite[Theorem 14.1.3 (ii)]{HilNeeb12}, $K_0 = K\cap G_0$ is connected
   and $K_0$ is a maximal compact subgroup of $G_0$. Moreover,
  $K$ intersects each connected component of $G$ and  $K\slash K_0 \cong G\slash G_0$.
    By the classical Lie theory, the Lie algebra of a compact Lie group is a direct product of
  an abelian Lie algebra and some simple Lie algebras.
  Then since the Lie algebra $\mathrm{Lie}(G)$ of $G$ is solvable and $K$ is
   compact, the Lie algebra $\mathrm{Lie}(K)  \subset \mathrm{Lie}(G)$ must be abelian.
   This implies that $K_0$ is a torus and hence a maximal torus in $G_0$.
  \nn

  In addition, since $G$ is virtually solvable, $G_0$ is actually
  solvable. This is because the radical $R$ of $G$ is a normal
  subgroup of $G_0$ and $\dim(R)=\dim(G_0)$ (since $G$ is virtually
  solvable). So $G_0\slash R$ is discrete. Then since $G_0$ is
  connected, $G_0$ must equal $R$.\nn

   Let $Z(G)$ be the center of $G$ and define $C = Z(G)\cap K$. Then $C$ is clearly a normal subgroup
    of $G$. Let $G'=G\slash C$ and $K' = K\slash C$ and let $\rho: G \rightarrow G'$ be the quotient map.
    Then since $\Gamma\cap K =\{ 1 \}$, $\Gamma\cong \rho(\Gamma) \subset
    G'$, we can think of $\Gamma$ as a subgroup of $G'$. So we
    have
  \begin{equation} \label{Equ:Equiv}
    M=\Gamma \backslash G \slash K \cong  \Gamma \backslash G' \slash K' .
  \end{equation}
  Let $G'_0=\rho(G_0)$ be the identity component of $G'$.
  Then $G'_0$ is a finite index normal subgroup of $G'$ and $G'_0$ is solvable.
    \begin{equation} \label{Equ:Quotient-0}
     \mathrm{Lie}(G'_0) = \mathrm{Lie}(G')= \mathrm{Lie}(G)\slash \mathrm{Lie}(C) =
    \mathrm{Lie}(G_0)\slash \mathrm{Lie}(C).
    \end{equation}
 Besides, let
  $K'_0 = K'\cap G'_0$ which is a maximal torus of
   $G'_0$ and we have
    \begin{equation} \label{Equ:Quotient-1}
     \mathrm{Lie}(K'_0) = \mathrm{Lie}(K') = \mathrm{Lie}(K) \slash \mathrm{Lie}(C)
      = \mathrm{Lie}(K_0)\slash \mathrm{Lie}(C) .
    \end{equation}
    \n

 \noindent\textbf{Claim-1:} $G'_0$ is linear and so $G'$ is linear.
 \nn

  A group is called \emph{linear} if it admits a faithful finite-dimensional representation.
   By~\cite[Theorem 16.2.9 (b)]{HilNeeb12}, a connected solvable Lie group $S$
   is linear if and only if $\mathfrak{t}\cap [\mathfrak{s}, \mathfrak{s}] =\{ 0 \}$ where
   $\mathfrak{s}$ and $\mathfrak{t}$ are Lie algebras of $S$ and its maximal torus $T_{S}$, respectively.
  And for a general connected solvable group $S$, the Lie subalgebra $\mathfrak{t}\cap [\mathfrak{s}, \mathfrak{s}]$ is
  always central in $\mathfrak{s}$.
  So for our $G_0$ and its maximal torus $K_0$,
  we have $\mathrm{Lie}(K_0)\cap [\mathrm{Lie}(G_0),
  \mathrm{Lie}(G_0)]$ is central in $\mathrm{Lie}(G_0)=\mathrm{Lie}(G)$.
  So $\mathrm{Lie}(K_0)\cap [\mathrm{Lie}(G_0),
  \mathrm{Lie}(G_0)] \subset \mathrm{Lie}(Z(G))$ and
  \begin{align} \label{Equ:Quotient-2}
    \mathrm{Lie}(K_0)\cap [\mathrm{Lie}(G_0),
  \mathrm{Lie}(G_0)]  \subset \mathrm{Lie}(K) \cap \mathrm{Lie}(Z(G))  =
   \mathrm{Lie}(C)
   \end{align}
 Then for the Lie group $G_0'$ and its maximal torus $K'_0$, we have
  \begin{align} \label{Equ:Quotient-3}
     \mathrm{Lie}(K'_0)\cap [\mathrm{Lie}(G_0'),
  \mathrm{Lie}(G_0')]= \frac{\mathrm{Lie}(K_0)}{\mathrm{Lie}(C)} \cap
  \left[\frac{\mathrm{Lie}(G_0)}{\mathrm{Lie}(C)}, \frac{\mathrm{Lie}(G_0)}{\mathrm{Lie}(C)}\right] =
  0. \quad
   \end{align}
  So by~\cite[Theorem 16.2.9 (b)]{HilNeeb12}, $G_0'$ is linear.
  Moreover, suppose $V$ is a faithful finite-dimensional representation of
  $G'_0$. Then $\R[G']\otimes_{\R[G'_0]} V$ is a
  faithful finite-dimensional representation of
  $G'$ where $\R[G']$ and $\R[G'_0]$ are the group rings of $G'$ and $G'_0$ over $\R$, respectively.
  So the Claim-1 is proved.\nn

 From the Claim-1 and~\eqref{Equ:Equiv}, we can just assume that our group $G$ is linear at the
  beginning. Under this assumption,
   $G_0$ is a connected linear solvable group.
  So there exists a simply connected solvable normal Lie subgroup $S$ of $G_0$ so that
  $G_0 = S\rtimes K_0$ and $[G_0,G_0]\subset S$ (see~\cite[Lemma 16.2.3]{HilNeeb12}).
   So
   $$\mathrm{Lie}(G_0) = \mathrm{Lie}(K_0) \oplus
   \mathrm{Lie}(S).$$
 More specifically, we can take $S = p^{-1}(V)$ where
 $p: G_0 \rightarrow G_0\slash [G_0,G_0]$ is the quotient map and
 $V$ is a vector subgroup of the abelian group $G_0\slash [G_0,G_0]$
 so that $G_0\slash [G_0,G_0] \cong  p(K_0) \times V$  (see the proof of~\cite[Theorem 16.2.3]{HilNeeb12}).
 Note that the vector subgroup $V$ is not unique, so $S$ is not unique either.
 \nn

  \noindent \textbf{Claim-2:} We can choose $S$ to be normal in $G$ and so $G \cong S\rtimes K$. \nn

  Indeed since $K$ is compact, we can choose a metric on $\mathrm{Lie}(G_0)$ which is invariant under
   the adjoint action of $K$. Then we can choose $V$ so that
   $\mathrm{Lie}(S)$ is orthogonal to $\mathrm{Lie}(K_0)$ in $\mathrm{Lie}(G_0)$.
  Then because $K_0$ is normal in $K$, the adjoint action of $K$ on $\mathrm{Lie}(G_0)$
   preserves $\mathrm{Lie}(K_0)$, so it also preserves the orthogonal complement $\mathrm{Lie}(S)$
   of $\mathrm{Lie}(K_0)$. This implies that $S$ is preserved under
  the adjoint action of $K$.\n

   Let $G_0$, $h_1 G_0, \cdots, h_m G_0$ be all the connected components of $G$.
   Since $K$ intersects each connected component of $G$, we can assume $h_i\in K$ for all
   $1\leq i \leq m$. Then any element $g \in G$ can be written as $g=g_0 h_i$
   for some $g_0\in G_0$ and $h_i\in K$. So
   $g S g^{-1} = g_0 h_i S h^{-1}_i g^{-1}_0 \subset g_0 S  g^{-1}_0 \subset
   S$. The Claim-2 is proved. \nn

  From the semidirect product $G=S\rtimes K$, we can define an injective group homomorphism
  $\alpha: G \rightarrow \mathrm{Aff}(S)= S\rtimes \mathrm{Aut}(S)$ as follows.
   For any $g\in G$, we can write $g = s_g k_g$ for a
    unique $s_g \in S$ and $k_g\in K$ since $S\cap K = S\cap K_0 = \{ 1\}$.
    Then $\alpha(g) : S \rightarrow S$ is the composition of the
  adjoint action of $k_g$ on $S$ and the left translation on $S$ by
  $s_g$, i.e. $\alpha(g) = L_{s_g}\circ \mathrm{Ad}_{k_g}$.  \nn

     \noindent \textbf{Claim-3:} $ \alpha(\Gamma) \backslash S$ is diffeomorphic to
   the double coset space $\Gamma\backslash G\slash K$.\nn

   Notice that each left coset in $G\slash K$ contains a unique element
   of $S$, so we have
    $$G\slash K = \{ s K \, ; \, s\in S\}.$$
  For any $\gamma\in \Gamma$, let $\gamma = s_{\gamma} k_{\gamma}$
   where $s_{\gamma} \in S$ and $k_{\gamma}\in K$, and we have
   \begin{equation}
       \gamma s K = s_{\gamma} k_{\gamma} s K = s_{\gamma} k_{\gamma} s k^{-1}_{\gamma} K =
     \alpha(\gamma) (s) K, \  \forall\, s\in S .
   \end{equation}
     So the natural action of $\Gamma$ on the left coset space $G\slash K$ can be identified with
      the action of $\alpha(\Gamma) \subset \mathrm{Aff}(S)$ on $S$.
      The Claim-3 is proved.
   \nn

    Let $\mathrm{Ad}: K \rightarrow \mathrm{Aut}(S)$
   denote the adjoint action of $K$ on $S$.
 Since $K$ is compact and $\mathrm{Ad}$ is continuous, so $\mathrm{Ad}(K) \subset \mathrm{Aut}(S) $ is also compact.
 Notice that $\alpha(\Gamma)$ is a subgroup of $S\rtimes  \mathrm{Ad}(K) \subset S\rtimes
 \mathrm{Aut}(S)$, so the closure $\overline{hol(\alpha(\Gamma))}$ of
 the holonomy group $hol(\alpha(\Gamma))$ in $\mathrm{Aut}(S)$ is contained in
 $\mathrm{Ad}(K)$. So $\overline{hol(\alpha(\Gamma))}$ is compact.
 This implies that $\Gamma \backslash G \slash K \cong \alpha(\Gamma) \backslash S$
 is an infra-solvmanifold in the sense of Def 2.
  \hfill $\qed$\nn

  \n

\noindent \textbf{Remark 1:}
  A simply-connected solvable Lie group is always linear, but for non-simply-connected
  solvable Lie groups, this is not always so.
  A counterexample is the quotient group of the Heisenberg group by
  an infinite cyclic group (see~\cite[p.169 Example 5.67]{HofMorris06}).
\\

\noindent $\textbf{\S 4.}$ \textbf{Two other definitions}\nn

 In addition to the three definitions listed at
the beginning, there are
 two other definitions of compact infra-solvmanifolds.\nn

 \begin{itemize}
   \item[Def 4:] A \emph{compact infra-solvmanifold} $M$ is an orbit space of the
  form $M=\Delta\backslash S$ where $S$ is a connected, simply connected solvable Lie group
   acted upon by a torsion-free cocompact
  closed subgroup $\Delta \subset \mathrm{Aff}(S)$ satisfying\n
     \begin{itemize}
        \item the identity component $\Delta_0$ of $\Delta$ is contained in the nil-radical of
        $S$,
        \n
        \item the closure of $hol(\Delta)$ in $\mathrm{Aut}(S)$ is
        compact.
   \end{itemize}
    \nn
   \item[Def 5:] A \emph{compact infra-solvmanifold} is a manifold of the form $\Delta \backslash G$,
   where $G$ is a connected, simply connected solvable Lie group, and $\Delta$
   is a torsion-free cocompact subgroup of $\mathrm{Aff}(G)$ so that the closure of $hol(\Delta)$ in $\mathrm{Aut}(G)$ is
        compact.
  \end{itemize}
\nn

 Note that the subgroup $\Delta$ in Def $4$ and Def $5$
  is not necessarily discrete. The Def $4$ is
  from~\cite[Definition 1.1]{FarJone97} and Def $5$ is from~\cite[Definition 1.1]{Baube04}. \nn

   It is clear that Def 2 $\Rightarrow$ Def $4$ $\Rightarrow$ Def $5$.
  Moreover, the results in~\cite{Baube04} implies that Def $5$ $\Rightarrow$ Def
  $2$.
  Indeed, for the manifold $\Delta\backslash G$ in Def
  $5$, its fundamental group is $\Gamma = \Delta\slash \Delta_0$ which is
   torision-free virtually poly-cyclic (see~\cite{FarJone90}).
  It is shown in~\cite{Baube04} that
  such a group $\Gamma$ determines a virtually solvable real linear algebraic group
   $H_{\Gamma}$ which contains $\Gamma$ as a discrete and Zariski-dense subgroup.
     $H_{\Gamma}$ is called the \emph{real algebraic hull} of $\Gamma$ in~\cite{Baube04}.
   In addition, we have $H_{\Gamma} = U \rtimes T$ where
   $T$ is a maximal reductive subgroup of $H_{\Gamma}$ and $U$ is the
   unipotent radical of $H_{\Gamma}$. The splitting of $H_{\Gamma}$ gives an injective group
   homomorphism $\alpha: H_{\Gamma} \rightarrow \mathrm{Aff}(U)$ and
   a corresponding affine action of $\Gamma < H_{\Gamma}$ on $U$. It is shown
   in~\cite{Baube04} that $M_{\Gamma} = \alpha(\Gamma) \backslash U$
    is a compact infra-solvmanifold whose fundamental group is $\Gamma$.
 $M_{\Gamma}$ is called the \emph{standard $\Gamma$-manifold}.
 Now since $\alpha(\Gamma)$ is a discrete subgroup of
 $\mathrm{Aff}(U)$, so $M_{\Gamma}$ is a compact infra-solvmanifold in the sense of Def $2$.
  Moreover by~\cite[Theorem 1.4]{Baube04}, $\Delta\backslash G$
 is diffeomorphic to $M_{\Gamma}$. \nn

 By the preceding arguments,
 the five definitions Def 1\,--\,Def 5 of compact infra-solvmanifolds
 are all equivalent.\\

  \noindent \textbf{Remark 2:}
 If we remove the ``cocompact'' in Def 1\,--\,Def 5, we may
 get noncompact infra-solvmanifolds, which are
 vector bundles over some compact infra-solvmanifolds
 (see~\cite[Theorem 6]{Wilk00}).\\

 \section*{Acknowledgements}
   The authors want to thank M.~Masuda, J.~B.~Lee, W.~Tuschmann and J.~Hillman for some
  helpful information and comments. Besides, the first author
  would like to thank Y. Karshon for providing him excellent circumstances to do
  research.\\

\end{document}